\documentclass[10pt,twoside,francais]{article}
\usepackage{babel}
\usepackage{graphicx}
\usepackage{amsmath,amsthm,amssymb,amsfonts}
\usepackage{Latex-document}
\def\thebibliography#1{\section*{R\'ef\'erences}
\list{[\arabic{enumi}]}{\settowidth \labelwidth{[#1]} \leftmargin
\labelwidth \advance \leftmargin \labelsep \usecounter{enumi}}
\def\newblock{\hskip .11em plus .33em minus .07em} \sloppy
\clubpenalty 4000 \widowpenalty 4000 \sfcode`\.=1000 \relax}

\def\ly{\fontencoding{U}\fontfamily{lasy}\fontseries{m}\fontshape{n}\selectfont}
\def\guil#1{\leavevmode\hbox{{\ly(\kern-0.20em(\kern+0.20em}}\nobreak{}\,#1\,%
  \nobreak\hbox{{\ly\kern+0.20em)\kern-0.20em)}}}
\def\thesection{\Alph{section}.}
\makeatletter
\gdef\thmhead@plain#1#2#3{%
  \thmname{#1}\thmnumber{\@ifnotempty{#1}{ }#2}%
  \thmnote{ {\mdseries#3}}}
\let\thmhead\thmhead@plain
\makeatother
\theoremstyle{plain}
\newtheorem{theoreme}{Th\'eor\`eme}
\newtheorem{proposition}[theoreme]{Proposition}

\newtheorem{corollaire}[theoreme]{Corollaire}
\theoremstyle{definition}
\newtheorem{definition}[theoreme]{D\'efinition}
\newtheorem*{exemple}{Exemple}
\def\C{\mathbf{C}}
\def\D{\mathbf{D}}

\def\R{\mathbf{R}}
\def\S{\mathbf{S}}
\def\T{\mathbf{T}}
\def\Z{\mathbf{Z}}
\def\id{\mathop{\mathrm{id}}}
\def\MT{\mathop{\Sigma}\nolimits}
\def\OB{\mathop{\overline{\Sigma}}\nolimits}

\def\Int{\mathop{\mathrm{Int}}\nolimits}

\def\wh{\hat}
\def\Alinea#1{\hfill\break%
  \hbox to \parindent{\hss{\textup{#1}}\enspace}\ignorespaces}
\def\alinea#1{\noindent%
  \hbox to \parindent{\hss{\textup{#1}}\enspace}\ignorespaces}
\def\up{\textup}
\def\op{\mathopen}

\def\from{\colon}
\def\res#1{\,\vert\,{}_{#1}}

\def\abs#1{\lvert #1 \rvert}

\def\Chi{\setbox0=\hbox{$\chi$} \mathord{\raise\dp0\hbox{$\chi$}}}

\def\db{\overline{\partial}}
\def\dc{\partial}

\title{\bf G\'eom\'etrie de Contact: \vskip -2mm de la Dimension Trois \vskip -2mm vers les Dimensions
Sup\'erieures\vskip 6mm}

\author{Emmanuel Giroux\vspace*{-0.5cm}\thanks{Unit\'e de Math\'ematiques Pures et Appliqu\'ees,
\'Ecole Normale Sup\'erieure de Lyon, 46 all\'ee d'Italie, 69364 Lyon cedex 07, France. M\'el:
giroux@umpa.ens-lyon.fr}}

\date{\vspace{-8mm}}

\begin{document}

\maketitle

\thispagestyle{first}

\setcounter{page}{405}

\begin{abstract}

\vskip 3mm

On d\'ecrit ici des relations entre la g\'eom\'etrie globale des vari\'et\'es de contact closes et celle de
certaines vari\'et\'es symplectiques, \`a savoir les vari\'et\'es de Stein compactes. L'origine de ces relations
est l'existence de livres ouverts adapt\'es aux structures de contact.

\vskip 4.5mm


\noindent {\bf Mots cl\'es:} Structures de contact, Livres
ouverts.
\end{abstract}

\vskip 12mm

La g\'eom\'etrie de contact en dimension trois a connu un essor important durant
la derni\`ere d\'ecennie gr\^ace au d\'eveloppement de m\'ethodes topologiques
ad\'equates. Dans le prolongement des travaux de D.~Bennequin~\cite{Be:eep} et
de Y.~Eliashberg~\cite{El:ocs}, la th\'eorie des \guil{surfaces convexes}~\cite
{Gi:ctc} et l'\'etude des rocades (bypasses)~\cite{Ho:cs1} ont men\'e \`a une
classification compl\`ete des structures de contact sur quelques vari\'et\'es
simples et, plus r\'ecemment, \`a une classification grossi\`ere sur toutes les
vari\'et\'es closes \cite{Co:isc,HKM,CGH}. En fait, comme on essaiera de le
montrer plus loin, les structures de contact en dimension trois sont des objets
purement topologiques, un peu comme les structures symplectiques en dimension
deux. En termes pr\'ecis, sur toute vari\'et\'e close~$V$ de dimension trois,
les classes d'isotopie des structures de contact se trouvent en correspondance
bijective avec les classes d'isotopie et de stabilisation des livres ouverts
dans~$V$, l'op\'eration \'el\'ementaire de stabilisation \'etant un plombage
positif~\cite{Gi:lo3}.

En dimension sup\'erieure, des m\'ethodes radicalement diff\'erentes permettent
de mettre en \'evidence une correspondance similaire~\cite{GM} et, au-del\`a,
de faire appara\^{\i}tre des liens \'etroits entre la g\'eom\'etrie globale des
vari\'et\'es de contact closes et celle de certaines vari\'et\'es symplectiques
compactes. Les livres ouverts qu'on associe \`a une structure de contact sont en
effet particuliers~: leurs pages sont des vari\'et\'es de Stein compactes, leur
monodromie est un diff\'eomorphisme symplectique \`a support dans l'int\'erieur
et l'op\'eration \'el\'ementaire de stabilisation qui les unifie est un plombage
lagrangien positif. En outre, l'outil essentiel pour les construire est la
th\'eorie des fibr\'es positifs que S.~Donaldson a introduite et d\'evelopp\'ee
en g\'eom\'etrie symplectique dans~\cite{Do:svs,Do:slp} et qui a \'et\'e
adapt\'ee en g\'eom\'etrie de contact dans~\cite{IMP}.

\section{Structures de contact et livres ouverts}

\vskip-5mm \hspace{5mm}

Dans ce texte, $V$ d\'esigne toujours une vari\'et\'e close et orient\'ee. Les
champs d'hyperplans tangents qu'on consid\`ere sur~$V$ sont coorient\'es, donc
aussi orient\'es puisque $V$ l'est. Un tel champ~$\xi$ est le noyau d'une forme
$\alpha$, appel\'ee \emph{\'equation} de~$\xi$, unique \`a multiplication pr\`es
par une fonction positive. On dit que $\xi$ est une \emph{structure de contact}
si $d\alpha$ induit sur~$\xi$ en tout point une forme symplectique directe,
\emph{i.e.} si $V$ est de dimension impaire $2n+1$ et si $\alpha \wedge
(d\alpha)^n$ est en tout point un \'el\'ement de volume direct pour
l'orientation de~$V$.

D'autre part, un \emph{livre ouvert} dans~$V$ est un couple $(K,\theta)$
form\'e des objets suivants~:
\begin{itemize}
\item
une sous-vari\'et\'e close $K \subset V$ de codimension deux \`a fibr\'e normal
trivial\,;
\item
une fibration $\theta \from V \setminus K \to \S^1$ qui, dans un voisinage $K
\times \D^2$ de $K = K \times \{ 0 \}$, co\"{\i}ncide avec la coordonn\'ee
angulaire normale.
\end{itemize}
On peut aussi voir les livres ouverts autrement. Soit $\phi \from F \to F$ un
diff\'eomorphisme d'une vari\'et\'e compacte \'egal \`a l'identit\'e pr\`es du
bord $K = \partial F$. Sa suspension, \`a savoir la vari\'et\'e compacte
$$ \MT(F,\phi) = \bigl( F \times [0,1] \bigr) \big/ \mathord\sim, \quad
   \text{o\`u \ $(p,1) \sim (\phi(p),0)$,} $$
est bord\'ee par $K \times \S^1$ \ --~car $\phi \res K = \id$~-- \ et la
vari\'et\'e close
$$ \OB(F,\phi) = \MT(F,\phi) \cup_\partial (K \times \D^2), $$
poss\`ede un livre ouvert \'evident. En outre, tout livre ouvert $(K,\theta)$
dans~$V$ identifie~$V$ \`a $\OB(F,\phi)$, o\`u $F$ est une fibre de~$\theta$ (un
peu r\'etr\'ecie) et $\phi$ l'application de premier retour sur~$F$ d'un flot
transversal aux fibres de~$\theta$ et constitu\'e, pr\`es de~$K$, de rotations
autour de~$K$. Le diff\'eomorphisme~$\phi$, d\'efini seulement \`a conjugaison
et isotopie pr\`es, est la \emph{monodromie} de $(K,\theta)$.

\smallskip

Toute la discussion \`a venir tourne autour de la d\'efinition suivante~:

\begin{definition}[\cite{Gi:lo3,GM}] \label{d:porteur} {\it
Une structure de contact~$\xi$ sur~$V$ est dite \emph{port\'ee} par un livre
ouvert $(K,\theta)$ si elle admet une \'equation~$\alpha$ ayant les
propri\'et\'es
suivantes~:
\begin{itemize}
\item
$\alpha$ induit sur~$K$ une forme de contact\,;
\item
$d\alpha$ induit sur chaque fibre~$F$ de~$\theta$ une forme symplectique\,;
\item
l'orientation de~$K$ d\'efinie par la forme de contact~$\alpha$ co\"{\i}ncide
avec son orientation comme bord de la vari\'et\'e symplectique $(F,d\alpha)$.
\end{itemize}
Une telle forme $\alpha$ sera dite \emph{adapt\'ee} \`a
$(K,\theta)$.}
\end{definition}

\begin{exemple}[\cite{GM}] \label{x:porteurs}
Soit $f \from (\C^n,0) \to (\C,0)$ une fonction holomorphe ayant \`a l'origine
un point critique isol\'e et soit $H$ l'hypersurface (singuli\`ere) $f^{-1}(0)$.
Il existe une boule ferm\'ee lisse~$B$ autour de l'origine dans~$\C^n$ et un
feuilletage de $B \setminus \{0\}$ par des sph\`eres strictement pseudoconvexes
$S_r$, o\`u $r \in \op]0,1]$ et $S_1 = \partial B$, tels que, pour $r$ assez
petit, les propri\'et\'es suivantes soient satisfaites~:
\begin{itemize}
\item
la sph\`ere $S_r$ est transversale \`a $H$, de sorte que $K = H \cap S_r$ est
une sous-vari\'et\'e close de~$S_r$ de codimension deux et \`a fibr\'e normal
trivial\,;
\item
l'application $\theta = \arg f \from S_r \setminus K \to \S^1$ est une fibration
qui fait de $(K,\theta)$ un livre ouvert\,;
\item
le livre ouvert $(K,\theta)$ porte la structure de contact sur~$S_r$ d\'efinie par le champ des tangentes complexes.
\end{itemize}
Autrement dit, chaque livre ouvert donn\'e dans la sph\`ere par le th\'eor\`eme
de fibration de J.~Milnor porte, \`a isotopie pr\`es, la structure de contact
standard.
\end{exemple}

\section{Structures de contact et livres ouverts en dimension trois}

\vskip-5mm \hspace{5mm}

En dimension trois, divers travaux ont depuis longtemps fait appara\^\i tre des
connivences entre les structures de contact et les livres ouverts sans toutefois
\'etablir aucun lien formel. Dans~\cite{TW}, W.~Thurston et H.~Winkelnkemper
construisent des formes de contact sur toute vari\'et\'e close~$V$ \`a partir
d'un livre ouvert dans~$V$. Avec les termes de la d\'efinition~\ref{d:porteur},
ils d\'emontrent en fait que tout livre ouvert dans~$V$ porte une structure de
contact. Dans~\cite{Be:eep} d'autre part, pour transformer en th\'eor\`eme de
g\'eom\'etrie de contact son r\'esultat sur les tresses ferm\'ees, D.~Bennequin
met en \'evidence la propri\'et\'e suivante~: toute courbe transversale \`a la
structure de contact standard~$\xi_0$ dans~$\R^3$ --~structure d'\'equation
$dz + r^2 d\theta = 0$~-- est isotope, parmi les courbes transversales, \`a une
tresse ferm\'ee c'est-\`a-dire une courbe transversale au livre ouvert form\'e
par l'axe des~$z$ et la coordonn\'ee angulaire~$\theta$. Or cette propri\'et\'e
vient de ce que ce livre ouvert porte~$\xi_0$. Enfin, dans~\cite{To:ccs},
I.~Torisu a clairement d\'egag\'e les relations entre les livres ouverts et les
configurations de th\'eorie de Morse consid\'er\'ees dans~\cite{Gi:ctc} pour
\'etudier les structures de contact \emph{convexes} au sens de~\cite{EG}.

\smallskip

La premi\`ere observation qui montre l'\'etroitesse des liens impos\'es par la
d\'efini\-tion \ref{d:porteur} et d\'ecoule de la stabilit\'e des structures de
contact est la suivante~:

\begin{proposition}[\cite{Gi:lo3}] \label{p:isotopie}
Sur une vari\'et\'e close de dimension trois, toutes les  structures de contact
port\'ees par un m\^eme livre ouvert sont isotopes.
\end{proposition}

Quant \`a la question de savoir quelles structures de contact poss\`edent un
livre ouvert porteur, la r\'eponse est simple~:

\begin{theoreme}[\cite{Gi:lo3}] \label{t:existence3}
Sur une vari\'et\'e close de dimension trois, toute structure de contact est
port\'ee par un livre ouvert.
\end{theoreme}

Cependant, comme l'illustre l'exemple des fibrations de Milnor, le livre ouvert
qui porte une structure de contact donn\'ee est loin d'\^etre  unique --~m\^eme
\`a isotopie pr\`es. Pour appr\'ehender ce ph\'enom\`ene, quelques d\'efinitions
sont utiles.

Soit $F \subset V$ une surface compacte \`a bord et $C \subset F$ un arc simple
et propre. On dit qu'une surface compacte $F' \subset V$ s'obtient \`a partir
de~$F$ par le \emph{plombage positif} (resp. \emph{n\'egatif}) d'un anneau le
long de $C$ si $F' = F \cup A$ o\`u $A \subset V$ est un anneau ayant les
propri\'et\'es suivantes~:
\begin{itemize}
\item
$A \cap F$ est un voisinage r\'egulier de~$C$ dans~$F$\,;
\item
$A$ est inclus dans une boule ferm\'ee~$B$ dont l'intersection avec~$F$ est
r\'eduite \`a $A \cap F$ et l'enlacement des deux composantes de $\partial A$
dans~$B$ vaut $1$ (resp.~$-1$).
\end{itemize}
Un r\'esultat de J.~Stallings affirme que, si $(K,\theta)$ est un livre ouvert
dans~$V$ et si $F$ est l'adh\'erence d'une fibre de~$\theta$, alors, pour toute
surface~$F'$ obtenue \`a partir de~$F$ par le plombage d'un anneau, il existe un
livre ouvert $(K',\theta')$ tel que $K'$ soit le bord de~$F'$ et que $F'$ soit
l'adh\'erence d'une fibre de~$\theta'$. Dans la suite, on dira que le livre
$(K',\theta')$ et l'entrelacs~$K'$ sont eux-m\^emes obtenus par plombage \`a
partir respectivement de $(K,\theta)$ et de~$K$. En outre, on dira qu'un livre
ouvert $(K',\theta')$ est une \emph{stabilisation} d'un autre $(K,\theta)$ s'il
s'obtient \`a partir de $(K,\theta)$ par une suite finie de plombages positifs.

\begin{theoreme}[\cite{Gi:lo3}] \label{t:unicite3}
Dans une vari\'et\'e close de dimension trois, deux livres ouverts quelconques
qui portent une m\^eme structure de contact ont des stabilisations isotopes.
\end{theoreme}

Les th\'eor\`emes \ref{t:existence3} et \ref{t:unicite3} permettent de traduire
nombre de questions sur les structures de contact en questions sur les livres
ouverts, autrement dit sur les diff\'eomorphismes des surfaces compactes \`a
bord. En ce sens, ce sont les analogues des th\'eor\`emes de S.~Donaldson~\cite
{Do:slp} sur les pinceaux de Lefschetz dans les vari\'et\'es symplectiques de
dimension quatre. Ils admettent cependant, \`a la diff\'erence de ceux-ci, des
d\'emonstrations purement topologiques dont on d\'ecrit bri\`evement les id\'ees
ci-dessous, apr\`es avoir introduit l'outil essentiel. On supposera le lecteur
familier avec certaines notions de g\'eom\'etrie de contact en dimension trois
(structures de contact vrill\'ees/tendues, invariant de Thurston-Bennequin des
courbes legendriennes, surfaces $\xi$-convexes).

\medskip

On appelle \emph{cellule poly\'edrale} dans~$V$ l'image d'un poly\`edre convexe
compact euclidien par un plongement topologique. Une telle cellule poss\`ede une
structure affine induite par son param\'etrage et son \emph{int\'erieur} est,
par d\'efinition, l'image de l'int\'erieur \guil{intrins\`eque} du poly\`edre,
c'est-\`a-dire de son int\'erieur topologique dans son enveloppe affine. Une
\emph{cellulation poly\'edrale} de $V$ d\'esigne ici  un recouvrement fini de
$V$ par des cellules poly\'edrales ayant les propri\'et\'es suivantes~:
\begin{itemize}
\item
les int\'erieurs des cellules forment une partition de $V$\,;
\item
le bord de chaque cellule~$D$ est une union de cellules~$D_j$ et les inclusions
$D_j \to D$ sont affines\,;
\item
les cellules de dimension deux (et moins) sont lisses, \emph{i.e.} sont les
images de plongements lisses.
\end{itemize}
Les cellulations poly\'edrales ont cet avantage sur les triangulations d'\^etre
tr\`es faciles \`a subdiviser~: toute subdivision d'un sous-complexe se prolonge
trivialement. En outre, elles jouent un r\^ole cl\'e dans la d\'emonstration du
th\'eor\`eme de Reidemeister-Singer donn\'ee dans~\cite{Si:trs}, d\'emonstration
qui sert de guide pour \'etablir le th\'eor\`eme~\ref{t:unicite3}.

\smallskip

\begin{proof}[\bf Esquisse de la d\'emonstration du
th\'eor\`eme~\ref{t:existence3}] Soit $\xi$ une structure de
contact sur $V$. On construit d'abord dans~$(V,\xi)$ une
\emph{cellulation de contact}, c'est-\`a-dire une cellulation
poly\'edrale $\Delta$ ayant les propri\'et\'es suivantes~:
\begin{itemize}
\item[1)]
chaque cellule de dimension~$1$ est un arc legendrien\up{\,;}
\item[2)]
chaque cellule de dimension~$2$ est $\xi$-convexe et l'invariant de
Thurston-Bennequin de son bord vaut~$-1$\up{\,;}
\item[3)]
chaque cellule de dimension~$3$ est contenue dans le domaine d'une carte de
Darboux.
\end{itemize}
On \'epaissit ensuite le $1$-squelette $L$ de~$\Delta$ en une surface compacte
$\wh F$ (presque) tangente \`a~$\xi$ le long de~$L$ et on choisit un voisinage
r\'egulier~$W$ de~$L$ assez petit pour que $F = \wh F \cap W$ soit une surface
proprement plong\'ee dans~$W$. Quitte \`a prendre~$W$ plus petit, $\xi$ admet
une \'equation~$\alpha$ v\'erifiant les conditions suivantes~:
\begin{itemize}
\item
$d\alpha$ induit sur~$F$ une forme d'aire\,;
\item
$\alpha$ est non singuli\`ere sur $K = \partial F$ et oriente~$K$ comme le bord
de $(F,d\alpha)$.
\end{itemize}
D'autre part, pour toute cellule $D$ de dimension~$2$, la propri\'et\'e~2) dit
que le bord de $D \cap (V \setminus \Int W)$ intersecte~$K$ en deux points (\`a
isotopie pr\`es). Il en r\'esulte qu'il existe une fibration $\theta \from V
\setminus K \to \S^1$ ayant $\Int F$ pour fibre. Quitte \`a rogner~$W$, on peut
supposer que $W$ est une union de fibres de~$\theta$ sur lesquelles $d\alpha$
induit une forme d'aire. Il reste \`a voir que $\xi$ est isotope, relativement
\`a~$W$, \`a une structure de contact port\'ee par $(K,\theta)$. Le point cl\'e
est que $\xi$ est tendue sur $W^* = V \setminus \Int W$ et que $\partial W^*$
est une surface $\xi$-convexe dont le d\'ecoupage est fourni par~$K$.
\end{proof}

\begin{proof}[\bf \'Etapes de la d\'emonstration du th\'eor\`eme
\ref{t:unicite3}] Soit $\Delta$ une cellulation de contact de
$(V,\xi)$. On dira ici qu'un livre ouvert porteur $(K,\theta)$ est
\emph{associ\'e} \`a $\Delta$ si, comme dans la d\'emonstration du
th\'eor\`eme~\ref{t:existence3}, l'une des fibres de~$\theta$
contient le $1$-squelette de~$\Delta$ et se r\'etracte dessus par
une isotopie de contact. En imitant~\cite{Si:trs}, on montre
d'abord que tout livre ouvert porteur admet une stabilisation
associ\'ee \`a une cellulation de contact. On se ram\`ene ainsi
\`a consid\'erer le cas de deux livres ouverts porteurs associ\'es
\`a des cellulations de contact $\Delta_0$ et~$\Delta_1$ en
position g\'en\'erale. D'apr\`es~\cite{Si:trs}, $\Delta_0$ et
$\Delta_1$ poss\`edent une subdivision commune $\Delta_2$ qui
s'obtient, \`a partir de~$\Delta_0$ comme de~$\Delta_1$, par des
bissections. On d\'eforme alors $\Delta_2$, relativement \`a
l'union des $1$-squelettes de $\Delta_0$ et~$\Delta_1$, en une
cellulation v\'erifiant les propri\'et\'es 1) et~3) des
cellulations de contact et ayant des $2$-cellules $\xi$-convexes.
Il suffit ensuite de subdiviser le $2$-squelette de~$\Delta_2$
pour obtenir une cellulation de contact~$\Delta$ et on montre pour
finir que le livre ouvert associ\'e \`a $\Delta$ est une
stabilisation de ceux associ\'es \`a $\Delta_0$ et \`a~$\Delta_1$.
\end{proof}

\medskip

On discute maintenant quelques corollaires des th\'eor\`emes \ref{t:existence3}
et~\ref{t:unicite3}.

On rappelle d'abord qu'un th\'eor\`eme de M.~Hilden et J.~Montesinos affirme que
toute vari\'et\'e close~$V$ de dimension trois est un rev\^etement \`a trois
feuillets de la sph\`ere~$\S^3$ simplement ramifi\'e au-dessus d'un entrelacs
(simplement signifie que le degr\'e local aux points de ramification dans~$V$
vaut~deux). On obtient le m\^eme r\'esultat pour les vari\'et\'es de contact
closes~:

\begin{corollaire}[\cite{Gi:lo3}] \label{c:revetements}
Toute vari\'et\'e de contact close de dimension trois est un rev\^ete\-ment \`a
trois feuillets de la sph\`ere de contact standard $(\S^3,\xi_0)$ simplement
ramifi\'e au-dessus d'un entrelacs transversal \`a $\xi_0$.
\end{corollaire}

Un autre corollaire concerne la dynamique des flots de Reeb. Un \emph{flot de
Reeb} sur une vari\'et\'e de contact est un flot qui pr\'eserve la structure de
contact tout en lui \'etant transversal et en pointant du c\^ot\'e positif.
Un exemple typique est le flot g\'eod\'esique sur le fibr\'e cotangent unitaire
d'une vari\'et\'e riemannienne. Les flots de Reeb d'une structure de contact
donn\'ee~$\xi$ sont en bijection avec les \'equations de~$\xi$~: \`a toute forme
$\alpha$ correspond l'unique champ de vecteurs $\nabla_\alpha$ qui engendre le
noyau de $d\alpha$ et sur lequel $\alpha$ vaut~$1$. En prenant une \'equation
de~$\xi$ adapt\'ee \`a un livre ouvert porteur, on obtient~:

\begin{corollaire}[\cite{Gi:lo3}]
Sur toute vari\'et\'e de contact close de dimension~trois, il existe un flot de
Reeb qui admet une section de Poincar\'e-Birkhoff, c'est-\`a-dire une surface
compacte qui rencontre toutes les orbites, dont l'int\'erieur est transversal au
flot et dont chaque composante du bord est une orbite p\'eriodique.
\end{corollaire}

En fait, il n'est pas exclu que tout flot de Reeb admette une telle section
\cite{HWZ} (ceci impliquerait la conjecture de Weinstein selon laquelle tout
flot de Reeb a une orbite p\'eriodique) mais c'est l\`a un probl\`eme de nature
diff\'erente, certainement inaccessible par des  m\'ethodes topologiques.

\smallskip

Une question naturelle au vu des th\'eor\`emes \ref{t:existence3} et~\ref
{t:unicite3} est de savoir comment lire sur la monodromie de ses livres ouverts
porteurs si une structure de contact est tendue, ou remplissable en un
quelconque sens. La seule r\'eponse connue concerne les structures de contact
\emph{holomorphiquement remplissables}, c'est-\`a-dire r\'ealisables comme
champs des tangentes complexes au bord de vari\'et\'es de Stein compactes. Le
corollaire suivant pr\'ecise un r\'esultat de A.~Loi et R.~Piergallini~:

\begin{corollaire}[\cite{LP,Gi:lo3}] \label{c:remplissage}
Une structure de contact sur une vari\'et\'e close de dimension trois est
holomorphiquement remplissable si et seulement si elle est port\'ee par un livre
ouvert dont la monodromie est un produit de twists de Dehn \`a droite.
\end{corollaire}

Pour finir, on donne un corollaire de pure th\'eorie des n\oe uds. On appelle
ici \emph{entrelacs fibr\'e} dans $V$ tout entrelacs orient\'e~$K$ pour lequel
il existe une fibration $\theta \from V \setminus K \to \S^1$ qui fait de $(K'
\theta)$ un livre ouvert et induit sur $K$ l'orientation prescrite. Lorsque $V$
est une sph\`ere d'homologie, un th\'eor\`eme de F.~Waldhausen assure que cette
fibration, si elle existe, est unique \`a isotopie pr\`es. Le r\'esultat suivant
r\'epond \`a une question pos\'ee par J.~Harer dans~\cite{Ha:cfl}~:

\begin{corollaire}[\cite{Gi:lo3}] \label{c:entrelacs}
Deux entrelacs fibr\'es quelconques dans une sph\`ere d'homologie enti\`ere~$V$
s'obtiennent l'un \`a partir de l'autre  par une suite de plombages et de \guil
{d\'eplombages} \up(op\'erations inverses\up).
\end{corollaire}

\begin{proof}[\bf D\'emonstration.]
Une trivialisation de~$V$ \'etant choisie, les classes d'homotopie de champs de
plans tangents \`a~$V$ sont rep\'er\'ees par leur \emph{invariant de Hopf}, \`a
savoir l'enlacement des fibres des applications $V \to \S^2$ correspondantes. On
consid\`ere alors un livre ouvert quelconque $(K,\theta)$ dans~$V$ et on note
$(K',\theta')$ un livre ouvert obtenu \`a partir de $(K,\theta)$ par un plombage
n\'egatif. Les trois observations suivantes d\'emontrent le corollaire~:
\begin{itemize}
\item
toute structure de contact~$\xi'$ port\'ee par $(K',\theta')$ est vrill\'ee car
l'\^ame de l'anneau plomb\'e est isotope \`a une courbe legendrienne non nou\'ee
dans ~$(V,\xi')$ dont l'invariant de Thurston-Bennequin vaut~$+1$\,;
\item
l'invariant de Hopf de $\xi'$ est sup\'erieur d'une unit\'e \`a celui des
structures de contact port\'ees par $(K,\theta)$ (voir~\cite{NR})\,;
\item
si deux structures de contact vrill\'ees ont le m\^eme invariant de Hopf, elles
sont isotopes d'apr\`es~\cite{El:ocs} et deux livres ouverts quelconques qui les
portent ont donc des stabilisations isotopes.
\end{itemize}
Cet argument borne en outre par $h+2$ le nombre des (d\'e)\,plombages n\'egatifs
n\'ecessaires pour passer d'un livre ouvert \`a un autre, o\`u $h$ d\'esigne la
diff\'erence entre les invariants de Hopf correspondants.
\end{proof}

\section{Structures de contact et livres ouverts en dimension sup\'erieure}

\vskip-5mm \hspace{5mm}

En dimension sup\'erieure \`a trois, les livres ouverts porteurs de structures
de contact ne sont pas quelconques~: leurs fibres ont une structure symplectique
invariante par la monodromie. Pour pr\'eciser ce point, quelques d\'efinitions
sont utiles.

Soit $F$ une vari\'et\'e compacte, \`a bord $K = \partial F$. Une forme
symplectique exacte~$\omega$ sur $\Int F$ est \emph{convexe \`a l'infini} s'il
existe sur $\Int F$ un champ de Liouville (champ de vecteurs $\omega$-dual d'une
primitive de~$\omega$) qui est transversal \`a toutes les hypersurfaces $K
\times \{t\}$, $t \in \op]0,1]$, o\`u $K \times [0,1]$ est un voisinage collier
de $K = K \times \{0\}$. On dit en outre que $(\Int F,\omega)$ est une \emph
{vari\'et\'e de Weinstein}~\cite{EG} s'il existe un tel champ de Liouville qui,
de plus, est le (pseudo)\,gradient d'une fonction de Morse $F \to \R$ constante
et sans points critiques sur $K$. L'exemple typique de vari\'et\'e de Weinstein
est l'int\'erieur d'une \emph{vari\'et\'e de Stein compacte}. On nomme ainsi
toute vari\'et\'e complexe compacte~$F$ qui admet une fonction strictement
pluri-sous-harmonique $f \from F \to \R$ constante et sans points critiques sur
le bord. La $2$-forme $i\dc\db f$ d\'efinit alors une structure symplectique. Il
ressort en fait du travail de Y.~Eliashberg~\cite{El:csm} que toute vari\'et\'e
de Weinstein est symplectiquement diff\'eomorphe \`a l'int\'erieur d'une telle
vari\'et\'e de Stein compacte.

Si maintenant $\alpha$ est une forme de contact adapt\'ee \`a un livre ouvert
$(K,\theta)$,  sa diff\'erentielle~$d\alpha$ induit sur chaque fibre de~$\theta$
une structure symplectique exacte convexe \`a l'infini. Celle-ci d\'epend du
choix de $\alpha$ mais sa compl\'etion~\cite{EG} est bien d\'efinie \`a isotopie
pr\`es. Le th\'eor\`eme de W.~Thurston et H.~Winkelnkemper et la proposition
\ref{p:isotopie} s'\'etendent alors ainsi en grande dimension~:

\begin{proposition}[\cite{GM}] \label{p:construction}
Soit $F$ une vari\'et\'e compacte avec, sur $\Int F$, une forme symplectique
exacte convexe \`a l'infini et soit $\phi \from F \to F$ un diff\'eomorphisme
symplectique \'egal \`a l'identit\'e pr\`es de~$K = \partial F$. Il existe
alors sur $\OB(F,\phi)$ une structure de contact port\'ee par le livre ouvert
\'evident. De plus, deux structures de contact port\'ees par un m\^eme livre
ouvert et qui induisent sur ses pages des structures symplectiques ayant des
compl\'etions isotopes sont isotopes.
\end{proposition}

Quant au th\'eor\`eme~\ref{t:existence3}, il se g\'en\'eralise comme suit~:

\begin{theoreme}[\cite{GM}] \label{t:existence}
Toute structure de contact sur une vari\'et\'e close~$V$  est port\'ee par un
livre ouvert dont chaque fibre est une vari\'et\'e de Weinstein.
\end{theoreme}

\begin{proof}[\bf Esquisse de la d\'emonstration]
Soit $\xi$ une structure de contact, $\alpha$ une \'equation de~$\xi$ et $J$ une
structure presque complexe sur $\xi$ calibr\'ee par $d\alpha \res \xi$. On note
$\nabla_\alpha$ le champ de Reeb associ\'e \`a~$\alpha$ et $g$ la m\'etrique
riemannienne sur $V$ qui vaut $d\alpha(.,J.)$ sur $\xi$ et rend $\nabla_\alpha$
unitaire et orthogonal \`a~$\xi$. En termes \'el\'ementaires, le th\'eor\`eme
principal de~\cite{IMP} montre qu'il existe des constantes $C, \eta > 0$ et des
fonctions $s_k \from V \to \C$, $k \ge 1$, v\'erifiant les conditions suivantes~:
\begin{itemize}
\item
en tout point de~$V$,
$$ \abs{ s_k(p) } \le C, \\\quad
   \abs{ ds_k - i k s_k \alpha } \le C k^{1/2} \quad \text{et} \quad
   \abs{ \db_\xi s_k } \le C\,; $$
\item
en tout point $p$ o\`u $\abs{ s_k(p) } \le \eta$,
$$ \abs{ \dc_\xi s_k (p) } \ge \eta k^{1/2} \,. $$
\end{itemize}
(Ici, $\dc_\xi^{} s_k$ et $\db_\xi^{} s_k$ sont les parties respectivement
$J$-lin\'eaire et $J$-anti\-li\-n\'eaire de $ds_k \res \xi$.) En termes plus
parlants, les fonctions $s_k$ sont des sections approximativement holomorphes
et \'equitransversales du fibr\'e $L^{\otimes k} \to V$, o\`u $L$ est le fibr\'e
hermitien trivial $V \times \C \to V$ muni de la connexion unitaire d\'efinie
par la forme $-i\alpha$.

Les estimations ci-dessus entra\^{\i}nent d'abord que, pour $\abs w  \le \eta$,
l'ensemble $K_w = s_k^{-1}(w)$ est une sous-vari\'et\'e et que la forme
$\alpha_w$ induite par~$\alpha$ sur $K_w$ est une forme de contact (voir~\cite
{IMP}). En effet, $\alpha_w$ est non singuli\`ere pour $k$ assez grand puisque
son noyau est \'egal au noyau de $ds_k \res \xi$ et que $\abs{ \dc_\xi^{} s_k }
\ge \eta k^{1/2}$ tandis que $\abs{ \db_\xi^{} s_k } \le C$. Mieux, ces
in\'egalit\'es montrent que, pour $k$ grand,le noyau de~$\alpha_w$ est proche
d'un sous-espace $J$-complexe de~$\xi$ si bien que $d\alpha_w$ y est non
d\'eg\'en\'er\'ee.

L'observation suivante est que l'application $\arg s_k \from V \setminus K \to
\S^1$ est une fibration dont les fibres sont transversales au champ de Reeb
$\nabla_\alpha$ en tout point o\`u $\abs{ s_k } \ge \eta$. Pour le voir, on
note que l'estimation sur $ds_k - i k s_k \alpha$ implique que
$$ \abs{ ds_k(\nabla_\alpha) - i k s_k } \le C k^{1/2} . $$
Ainsi, en un point $p$ o\`u $\abs{ s_k(p) } \ge \eta$ et pour $k$ assez grand,
$ds_k(p)\,(\nabla_\alpha)$ est proche de $iks_k(p)$, \emph{i.e.} est non nul et
presque orthogonal \`a $s_k(p)$. Par suite, les sous-vari\'et\'es
$$ s_k^{-1} (R_\theta), \quad \text{o\`u \ \ }
   R_\theta = \{ re^{i\theta}, \; r > \eta \}, $$
sont transversales au champ de Reeb $\nabla_\alpha$.

Ces arguments montrent que le livre ouvert $(K = K_0, \; \theta = \arg s_k)$,
pour $k$ assez grand, porte la structure de contact $\xi = \ker \alpha$. Il
reste \`a v\'erifier que les fibres de~$\theta$ sont des vari\'et\'es de
Weinstein. Pour simplifier, on prouve ci-dessous l'assertion analogue en
g\'eom\'etrie symplectique.
\end{proof}

\begin{proposition} \label{p:weinstein}
Soit $W$ une vari\'et\'e close, $\omega$ une forme symplectique enti\`ere
sur~$W$ et $H_k$ une sous-vari\'et\'e symplectique de~$W$ en dualit\'e de
Poincar\'e avec $k\omega$ et obtenue par la construction de
Donaldson~\up{\cite{Do:svs}},
\`a partir d'un fibr\'e hermitien en droites $L$ muni d'une connexion unitaire
de courbure $-i\omega$. Pour $k$ assez grand, $(W \setminus H_k, \omega)$ est
une vari\'et\'e de Weinstein.
\end{proposition}

\begin{proof}[\bf D\'emonstration.]
En reprenant les arguments de~\cite{Do:slp}, on peut supposer que $H_k$ est le
lieu d'annulation d'une section $s_k \from V \to L^{\otimes k}$ qui v\'erifie,
en tout point de~$W$,
$$ \abs{ \db_k s_k } \le c \abs{ \dc_k s_k }  \quad \text{avec} \quad
   c < \frac 1 {\sqrt2} . $$
Dans la trivialisation de $L^{\otimes k}$ donn\'ee au-dessus de $W \setminus
H_k$ par la
section unitaire $u = s_k / \abs{ s_k}$, la connexion est d\'efinie par une
$1$-forme $-i\lambda$ o\`u $d\lambda = k\omega$. Si on pose $s_k = \varphi u$,
l'in\'egalit\'e ci-dessus donne
$$ \abs{ d\varphi / \varphi + J^*\lambda }
   < \abs{ d\varphi / \varphi - J^*\lambda }, $$
ce qui montre que $J^*\lambda$ est plus loin de $d\varphi / \varphi$ que de
$-d\varphi/\varphi$. Le champ de Liouville dual de $\lambda$ est alors un
pseudogradient de $\log\varphi$.
\end{proof}

\smallskip

Comme en dimension trois, le livre ouvert porteur d'une structure de contact
donn\'ee n'est pas unique. On d\'ecrit dans~\cite{GM} une op\'eration de \emph
{plombage le long d'un disque lagrangien} --~dans laquelle les twists de
Dehn-Seidel viennent remplacer les twists de Dehn~-- qui permet d'\'etablir des
analogues du th\'eor\`eme~\ref{t:unicite3} et du corollaire~\ref{c:remplissage}.
Ces r\'esultats ram\`enent l'\'etude des structures de contact \`a celles des
diff\'eomorphismes symplectiques des vari\'et\'es de Stein compactes qui sont
l'identit\'e pr\`es du bord. Ils permettent peut-\^etre ainsi de rapprocher les
travaux de Y.~Eliashberg, H.~Hofer et A.~Givental sur la th\'eorie symplectique
des champs de ceux de, par exemple, de P.~Seidel sur l'homologie de Floer et les
groupes de diff\'eomorphismes symplectiques. On peut aussi se demander si le
th\'eor\`eme~\ref{t:existence} cache des obstructions \`a l'existence d'une
structure de contact sur les vari\'et\'es closes. D'apr\`es~\cite{Qu:obd}, toute
vari\'et\'e close $V$ de dimension $2n+1$ poss\`ede un livre ouvert dont chaque
fibre a le type d'homotopie d'un complexe cellulaire de dimension~$n$. Il est
probable que, si $V$ admet un champ d'hyperplans tangents muni d'une structure
presque complexe, il existe un tel livre ouvert pour lequel chaque fibre est une
vari\'et\'e presque complexe et est donc, d'apr\`es~\cite{El:csm}, l'int\'erieur
d'une vari\'et\'e de Stein compacte. Toute la difficult\'e serait donc vraiment
de r\'ealiser la monodromie par un diff\'eomorphisme symplectique... Dans cet
ordre d'id\'ee, voici un corollaire concret du th\'eor\`eme~\ref{t:existence}
obtenu par F.~Bourgeois et qui montre, en r\'eponse \`a une vieille question,
que tout tore de dimension impaire poss\`ede une structure de contact~:

\begin{corollaire}[\cite{Bo:tcm}]
Si une vari\'et\'e close $V$ admet une structure de contact, $V \times \T^2$ en
admet une aussi.
\end{corollaire}

\begin{proof}[\bf D\'emonstration.]
Soit $\xi$ une structure de contact sur~$V$, soit  $\alpha$ une \'equation de
$\xi$ adapt\'ee \`a un livre ouvert porteur $(K,\theta)$ et soit $N = K \times
\D^2$ un voisinage de $K = K \times \{0\}$ dans lequel $\theta$ est la
coordonn\'ee angulaire normale. On note $r$ la coordonn\'ee radiale normale
dans~$N$ et on pose
$$ \tilde\alpha = \alpha + f(r) ( \cos\theta \, dx_1 - sin\theta \, dx_2 ),
   \qquad (x_1, x_2) \in \T^2 = \R^2\!/\Z^2, $$
o\`u la fonction$f(r)$ vaut $r$ pour $r \le r_0$, $1$ pour $r \ge 2r_0$ et
v\'erifie $f'(r) \ge 0$. Un calcul montre que, si on choisit $r_0$ assez petit,
$\tilde\alpha$ est une forme de contact sur $V \times \T^2$.
\end{proof}

\label{lastpage}

\end{document}